\documentclass[12pt]{article}
\usepackage{amsfonts}
\usepackage{amsmath}
\usepackage{amsthm}
\usepackage{subfigure}
\usepackage{graphics}
\numberwithin{equation}{section}

\newcommand{\N}{\mathbb{N}}
\newcommand{\wtil}{\widetilde}
\newcommand{\ds}{\displaystyle}
\newcommand{\Kr}{Krawtchouk\ }
\newcommand{\pf}{P$_{\textup V}$}
\newtheorem{theorem}{Theorem}[section]

\begin{document}

\title{The generalized Krawtchouk  polynomials and the fifth Painlev\'e equation}
\author{Lies Boelen, Galina Filipuk, Christophe Smet,  Walter Van Assche, \\ Lun Zhang}
\date{\today}
\maketitle

\begin{abstract}
We study the recurrence coefficients of the orthogonal polynomials
with respect to a semi-classical extension of the Krawtchouk weight.
We derive a coupled discrete system for these coefficients and show
that they satisfy the fifth Painlev\'e equation when viewed as
functions of one of the parameters in the weight.

\end{abstract}



\section{Introduction}

\subsection{Orthogonal polynomials}
It is well-known (see, for instance, \cite{Chihara, Ismail}) that a
sequence of orthonormal polynomials $(p_n)_{n \in \N}$, i.e.,
\begin{equation}\label{eq: orthonormality}
\int p_m(x)p_n(x)\,d\mu(x)=\delta_{m,n},
\end{equation}
where  $\mu$ is a positive measure with support
 on the real line and $\delta_{m,n}$ is the Kronecker delta, satisfies a three-term recurrence relation
\begin{equation}\label{3 term orthonormal}
 x p_n(x)=a_{n+1}p_{n+1}(x)+b_np_n(x)+a_np_{n-1}(x)
\end{equation}
with the recurrence coefficients given by the   following integrals
\begin{gather}\label{recu:orthonormal}
a_n=\int x p_n(x)p_{n-1}(x)\,d\mu(x),\qquad b_n=\int
xp_n^2(x)\,d\mu(x).
\end{gather}
Here it is assumed that $p_{-1}=0$.

We will consider discrete orthogonal polynomials for which the
weight is supported on an equidistant lattice $\{h n+n_0|\,n\in
A\subset \mathbb{Z}\}$ with parameters $h$ and $n_0$.  In case of
the lattice $\N_0=\{0,1,\ldots\}$ the orthogonality condition
(\ref{eq: orthonormality}) reads
\begin{equation}\label{eq: discrete orthonormality}
\sum_{k=0}^{\infty} p_m(k)p_n(k)\,w(k)=\delta_{m,n}.
\end{equation}
The examples of (classical) discrete orthogonal polynomials on an
equidistant lattice include the Meixner and Charlier polynomials
with the lattice $\N_0$ and the \Kr and Hahn polynomials ($N+1$
polynomials orthogonal on $\{0,1,$ $\ldots,$ $N\}$). For more
information about discrete (classical) orthogonal polynomials, we
refer to \cite{Ismail,NSU} (see also \cite{Meixner1, Charlier,
Meixner2} for the definitions of classical weights and their
semi-classical extensions and further references).

A useful characterization of classical polynomials is the Pearson
equation
$$ [\sigma(x)w(x)]'=\tau(x)w(x), $$ where $\sigma $ and $\tau$ are polynomials
satisfying $\deg\sigma\leq 2$ and $\deg\tau= 1$. In case of discrete
polynomials on a linear lattice the Pearson equation reads
$$ \nabla[\sigma(x)w(x)]=\tau(x)w(x), $$
where $\nabla$ is the backward difference operator
$$\nabla f(x)=f(x)-f(x-1).$$
 Semi-classical orthogonal polynomials are defined
as orthogonal polynomials for which the weight satisfies a Pearson
equation for which $\deg \sigma > 2$ or $\deg \tau \neq 1$ (see
\cite{HvR,Mar}).  The recurrence coefficients  in the three-term
recurrence relation (\ref{3 term orthonormal}) for classical
orthogonal polynomials (e.g., Hermite, Laguerre, Jacobi, Charlier,
Meixner, \Kr and others) can be found explicitly in contrast to
non-classical weights. The recurrence coefficients of semi-classical
weights usually obey nonlinear recurrence relations, which, in many
cases, can be identified as discrete Painlev\'e equations; see for
instance \cite{BV, Magnus} and the references therein. Moreover,
when viewed as functions of one of the parameters in the weight, the
recurrence coefficients satisfy  (continuous) Painlev\'e equations
(see \cite{Meixner1, Charlier, Meixner2} in case of the discrete
semi-classical Meixner and Charlier weights  and \cite{Laguerre} and
the references therein for other semi-classical continuous weights).
This paper gives another example of the connection between discrete
orthogonal polynomials and the Painlev\'e equations.

\subsection{Statement of the results} The classical \Kr polynomials
for each $N\in\N $ and $n \leq N$ are defined by
\begin{equation}
K_n(x)={}_2F_1\left(-n,-x;-N;\frac{1}{p}\right),\quad p\in(0,1),
\end{equation}
where
\begin{equation*}
{}_2F_1\left(a,b;c;z\right)=\sum_{s=0}^{\infty}\frac{(a)_s
(b)_s}{(c)_s s!}z^s ,
\end{equation*}
with the Pochhammer symbol $(a)_n$ being defined by
$$   (a)_n=\prod_{j=0}^{n-1}(a+j) = a (a+1)\cdots (a+n-1), $$
is the hypergeometric function; cf. \cite[Chapter 15]{DLMF}. They
satisfy orthogonality conditions with respect to the binomial weight
$$w(k)= \binom{N}{k}  p^k(1-p)^{N-k},\quad k\in\{ 0,1,\ldots,N \}.
$$
The polynomials satisfy (\ref{3 term orthonormal}) for
$n=0,\ldots,N-1$ with $K_{-1}=0$ and $a_0=0$. The \Kr polynomials
are, in fact, the Meixner polynomials
$$M_n(x)={}_2F_1\left(-n,-x;\beta;1-\frac{1}{\ell}\right)$$ with
$\beta=-N$ and $\ell=p/(p-1)$. The recurrence coefficients for the
orthonormal \Kr polynomials are given by
\begin{equation}  \label{KrawRec}
   a_n^2=np(1-p)(N+1-n),\quad b_n=p(N-n)+n(1-p).
\end{equation}
Observe that $a_0^2=0=a_{N+1}^2$.

In this paper, we are concerned with a semi-classical generalization
of Krawtchouk polynomials with the weight function defined by
\begin{equation}\label{potential explicit}
 w(x) = \binom{N}{x}\frac{c^x }{(1-\alpha)_x},\qquad
 x\in\{0,\ldots,N\},
\end{equation}
where $\alpha<1$ and $c> 0$ are two real parameters. Note that we
can recover the Krawtchouk weight from \eqref{potential explicit} by
letting $c$ and $-\alpha$ tend to infinity and  $-\frac{c}{\alpha}
\to \frac{p}{1-p}$. Our first result gives a discrete system
satisfied by the recurrence coefficients.

\begin{theorem}\label{Th1}
Let $a_n$ and $b_n$ be the recurrence coefficients in \eqref{3 term
orthonormal} for the weight \eqref{potential explicit}. We then have
that
\begin{align}
x_n &=\frac{1}{N}\left(\frac{a_n^2}{c}+n\right), \label{def:x_n}\\
y_n &=-\frac{b_n+N+1+c-n-\alpha}{N}, \label{def:y_n}
\end{align}
satisfy the following discrete system
\begin{equation}\label{gen system}\displaystyle
\begin{cases}
(x_n+y_n)(x_{n+1}+y_n)=\displaystyle-\frac{y_n(N+1+Ny_n)
(N+1-\alpha+Ny_n)}{c N},\\
\\
(x_n+y_n)(x_n+y_{n-1})=\displaystyle\frac{x_n(-N-1+Nx_n)
(\alpha-N-1+Nx_n) }{N(Nx_n-n )},
\end{cases}
\end{equation}
with initial conditions
\begin{align}\label{eq:initial value xy}
x_0=0,\quad
y_0=-\frac{N+1+c-\alpha}{N}-\frac{c}{1-\alpha}\frac{M(-N+1,2-\alpha,-c)}{M(-N,1-\alpha,-c)},
\end{align}
where $M(a, b, z)$ is the confluent hypergeometric function
$_1F_1(a; b; z)$ defined by
\begin{equation*}
_1F_1(a; b; z)=\sum_{s=0}^{\infty}\frac{(a)_s }{(b)_s s!}z^s;
\end{equation*}
(cf. \cite[Chapter 13]{DLMF}).
\end{theorem}

Recall that a truncated confluent hypergeometric function is a
Laguerre polynomial, in particular $M(-N,1-\alpha,-c) =
\frac{N!}{(1-\alpha)_N}L_N^{(-\alpha)}(-c)$; see \cite[Eq.
13.6.19]{DLMF}. Hence the initial value $y_0$ is in terms of a ratio
of Laguerre polynomials.

The system (\ref{gen system}) can be obtained by a limiting
procedure from $\alpha$-dP{$_{\textup IV}$} \cite{GR,Walter} given
by
\begin{equation}\label{dP4}
\begin{cases}
(X_n+Y_n)(X_{n+1}+Y_n)=\displaystyle\frac{(Y_n-\wtil A)(Y_n-\wtil
B)(Y_n-\wtil C)(Y_n-\wtil D)}
{(Y_n+\Gamma-Z_n)(Y_n-\Gamma-Z_n)},\\
\\
(X_n+Y_n)(X_n+Y_{n-1})=\displaystyle\frac{(X_n+\wtil A)(X_n+\wtil
B)(X_n+\wtil C)(X_n+\wtil
D)}{(X_n+\Delta-Z_{n+1/2})(X_n-\Delta-Z_{n+1/2})},
\end{cases}
\end{equation}
with $\wtil A+\wtil B+ \wtil C+\wtil D=0$. Indeed, by taking
$$   X_n=x_n-1/\varepsilon,\quad Y_n=y_n+1/\varepsilon, $$
$$ \wtil A=1/\varepsilon,\quad \wtil B=-3/\varepsilon+(2N+2-\alpha)/N,$$
$$ \wtil C=1/\varepsilon-
(N+1)/N, \quad \wtil D=1/\varepsilon-(N+1-\alpha)/N,$$
$$ Z_n=(2n-1)/(2N)+1/\varepsilon, \quad \Gamma^2=4c/(N\varepsilon),
\quad \Delta=2/\varepsilon,$$ and letting $\varepsilon$ tend to
zero, we get (\ref{gen system}).

It is worthwhile to point out that if we take $\alpha=0$ and
$c=\frac{p}{1-p}$, the weight function \eqref{potential explicit}
reduces to the one considered in \cite{LiesThesis}, where the author
derived a discrete system for the recurrence coefficients. We also
observe that the system in \cite{Smet}, which was also used in
\cite{Meixner2} for the generalized Meixner weight, can be obtained
from (\ref{gen system}) with an appropriate choice of the parameters
and scaling of $x_n$ and $y_n$.

Since there are two parameters in the weight \eqref{potential
explicit}, the recurrence coefficients are dependent on these
parameters as well. Our next theorem shows that, when viewed as a
function of $c$, the recurrence coefficients $a_n(c)$ and $b_n(c)$
are related to the fifth Painlev\'e equation.
\begin{theorem}\label{thm: P5}
With $x_n(c) $ and $y_n(c)$ defined in \eqref{def:x_n} and
\eqref{def:y_n}, we have
\begin{equation}\label{transformation xn}
x_n(c)=\frac{(cy'+n-N)^2-2(n-N)(n-N+c y')y-A_1
y^4+A_2y^3-A_3y^2}{4cN(y-1)y^2},
\end{equation}
\begin{equation}\label{transformation}
y_n(c)=\frac{N-n+y(1+n+c-\alpha+(\alpha-N-1)y)-cy' }{2N(y-1)y},
\end{equation}
where
$$A_1=(1+N-\alpha)^2,\;\;A_2=2(1+N-\alpha)(1+c+N-\alpha),\;\;A_3=(1+c+n-\alpha)(1+c-n+2N-\alpha)$$
and
 $y=y(c)$ is the solution of the fifth Painlev\'e equation \pf
\
\begin{equation}\label{P5}
y''  =
\left(\displaystyle\frac{1}{2y}+\displaystyle\frac{1}{y-1}\right)(y')^2
- \displaystyle\frac{y'}{c}+\displaystyle\frac{(y-1)^2}{c^2}\left(A
y + \displaystyle\frac{B}{y}\right)+\displaystyle\frac{C
y}{c}+\,\,\displaystyle\frac{D y (y+1)}{y-1},
\end{equation}
with the parameters given by
\begin{equation}\label{par}
A=\frac{(\alpha-N-1)^2}{2},\;\;B=-\frac{(n-N)^2}{2},\;\;
C=-(n+\alpha),\;\;D=-\frac{1}{2}.
\end{equation}
\end{theorem}
Theorem \ref{thm: P5} gives the relationship between $x_n,\,y_n$
(and hence $a_n^2,\; b_n$) and solutions of  \pf\; explicitly.

The Painlev\'e equations possess the so-called Painlev\'e property:
the solutions have no movable branch points. They were discovered by
Painlev\'{e} and his colleagues at the beginning of the twentieth
century while classifying all second-order ordinary differential
equations of the form
\begin{equation}\label{eq:2nd ODE}
w''=\mathcal{R}(z,w,w'),
\end{equation}
where $'=d/dz$, the function $\mathcal{R}$ is rational in $w$ and
$w'$, meromorphic in $z$, which possess the Painlev\'e property. It
turns out that, up to M\"{o}bius transformations, only fifty
equations of the form \eqref{eq:2nd ODE} have the Painlev\'e
property \cite{Gambier,Painleve1,Painleve2}. Forty-four of these
equations can either be linearized, be transformed to a Riccati
equation or be solved in terms of elliptic functions. The six
remaining equations are now known as the Painlev\'{e} equations,
which are often referred to as nonlinear special functions
\cite{Clarkson06} and have numerous applications in mathematics and
mathematical physics.

The Painlev\'{e} equations cannot be solved in terms of elementary
functions or known classical special functions in general. For
certain combinations of parameters, however, P$_{\textup
{II}}$--P$_{\textup {VI}}$ have solutions expressed in terms of
special functions. For \pf, the choice of parameters \eqref{par} is
exactly when P$_{\textup {V}}$ admits classical solutions
expressible in terms of confluent hypergeometric functions
(equivalently, Kummer functions or Whittaker functions), see
\cite[$\S$ 32.10 (v)]{DLMF}. In fact, it is the case when confluent
hypergeometric functions have the associated Laguerre polynomials as
special cases. Hence, the condition
 \eqref{par} on the parameters in P$_{\textup {V}}$ is actually that when P$_{\textup {V}}$ has rational solutions (see \cite[$\S$ 32.8 (v)]{DLMF}).
 This is also consistent with the
initial condition \eqref{eq:initial value xy} for $y_0$. For more
information about classical and rational solutions about P$_{\textup
{V}}$  and the associated $\tau$-functions, we refer to
\cite{Clarkson05,Clarkson06,Forrester,Gromak76,GrLSh,KLM,Mas04,MOK,NY,Oka,Watanabe}.

The rest of this paper is devoted to the proofs of Theorems
\ref{Th1} and \ref{thm: P5}. They are given in Sections
\ref{sec:proof of thm1} and \ref{sec:proof of p5}, respectively. In
Section~\ref{num} we study  recurrence coefficients numerically. We
conclude this paper with a discussion in Section \ref{sec:dis}.

\section{Proof of Theorem \ref{Th1}}\label{sec:proof of thm1}
The proof of Theorem \ref{Th1} relies on the ladder operators for
discrete orthogonal polynomials, which were studied in
\cite{Ismail_ladder}. We start with a brief description of this
aspect.

\subsection{Ladder operators for discrete orthogonal polynomials}

Given a weight function $w$, we define a potential
\begin{equation}\label{potential}
u(x)=-\frac{\nabla w(x)}{w(x)}=\frac{w(x-1)-w(x)}{w(x)},
\end{equation}
which is a discrete analogue of the external field generated by $w$.
The action of the forward difference operator
$$\Delta f(x)=f(x+1)-f(x)$$ on orthogonal polynomials $p_n $ is
given by
\begin{equation}\label{ladder p_n}
\Delta p_n(x)=A_n(x) p_{n-1}(x)-B_n(x)p_n(x).
\end{equation}
In case of a weight $w$ on the lattice $\N_0$ with $w(-1)=0$,  the
coefficients $A_n(x)$ and $B_n(x)$ in (\ref{ladder p_n}) are given
by
\begin{equation}\label{eq:An}
A_n(x)=a_n\sum_{k\in\N_0}
p_n(k)p_n(k-1)\frac{u(x+1)-u(k)}{x+1-k}w(k) ,
\end{equation}
and
\begin{equation}\label{eq:Bn}
B_n(x)=a_n\sum_{k\in\N_0}p_n(k)p_{n-1}(k-1)\frac{u(x+1)-u(k)}{x+1-k}w(k).
\end{equation}
In the case of a weight $w$ supported on a finite lattice
$\{0,1,\ldots,N\}$ with boundary conditions $w(-1)=0$ and
$w(N+1)=0$, the coefficients $A_n(x)$ and $B_n(x)$ in (\ref{ladder
p_n}) are given by
\begin{equation}\label{eq:An finite}
A_n(x)=a_n\frac{p_n(N+1)p_n(N)}{N-x}w(N)+ a_n\sum_{k=0}^{N}
p_n(k)p_n(k-1)\frac{u(x+1)-u(k)}{x+1-k}w(k) ,
\end{equation}
and
\begin{equation}\label{eq:Bn finite}
B_n(x)=a_n\frac{p_n(N+1)p_{n-1}(N)}{N-x}w(N)+
a_n\sum_{k=0}^{N}p_n(k)p_{n-1}(k-1)\frac{u(x+1)-u(k)}{x+1-k}w(k).
\end{equation}
In both cases, when the lattices are finite or infinite, the
following compatibility relations between the functions $A_n$ and
$B_n$ hold:
\begin{equation}\label{ladder eq1}
B_n(x)+B_{n+1}(x)=\frac{x-b_n}{a_n}A_n(x)-u(x+1)+\sum_{j=0}^n
\frac{A_j(x)}{a_j},
\end{equation}
\begin{equation}\label{ladder eq2}
a_{n+1}A_{n+1}(x)-a_n^2\frac{A_{n-1}(x)}{a_{n-1}}=(x-b_n)B_{n+1}(x)-(x+1-b_n)B_{n}(x)+1.
\end{equation}

\subsection{Proof of Theorem \ref{Th1}}\label{sec21}

To prove Theorem \ref{Th1}, we follow the theme in
\cite{LiesThesis}, which made use of the compatibility relations
between the ladder operators.


It is readily seen that, with $w$ defined in \eqref{potential
explicit}, $w(-1)=0$ and $w(N+1)=0$. By \eqref{potential}, we have
\begin{equation}\label{eq:potential}
u(x)=-1+\frac{x(x-\alpha)}{c(N+1-x)}.
\end{equation} It then follows
from \eqref{eq:An finite} and \eqref{eq:Bn finite} that
$$A_n(x)=a_nR_n\frac{x}{N-x}+a_nT_n\frac{1}{N-x},$$
where
$$R_n=\frac{1}{c}\sum_{k=0}^N p_n(k)p_n(k-1)w(k),$$
$$T_n=p_n(N+1)p_n(N)w(N)+\frac{1}{c}\sum_{k=0}^Np_n(k)p_n(k-1)\frac{Nk+(N+1)(1-\alpha)
}{N-k+1}w(k),$$ and
$$B_{n}(x)=\frac{x}{N-x}r_n+\frac{1}{N-x}t_n$$ with $t_0=0$ and
$$r_n=\frac{a_n}{c}\sum_{k=0}^{N}p_n(k)p_{n-1}(k-1)w(k),$$
$$t_n=a_np_n(N+1)p_{n-1}(N)w(N)+\frac{a_n}{c}\sum_{k=0}^Np_n(k)p_{n-1}(k-1)
\frac{Nk+(N+1)(1-\alpha)}{N+1-k}w(k).$$ Using the orthonormality  of
the polynomials one immediately gets that
$$R_n=\frac{1}{c},\qquad r_n=0,$$ since $p_n(k-1)=p_n(k)+$lower order
terms. Thus, we obtain
$$A_n(x)=a_n\frac{x}{c(N-x)}+a_n T_n\frac{1}{N-x},$$
$$B_n(x)=\frac{t_n}{N-x},$$ with $t_0=0.$

%
%

Next we use the compatibility relations (\ref{ladder eq1}) and
(\ref{ladder eq2}). The first compatibility relation (\ref{ladder
eq1}) gives rise to the following two equations, after collecting
coefficients of equal powers of $x$:
\begin{equation}\label{4.2.2}
c(t_n+t_{n+1}+b_n T_n)=-1+\alpha+cN+c\sum_{j=0}^{n}T_j,
\end{equation}
\begin{equation}\label{4.2.3}
c T_n-b_n-1+\alpha-c+n=0.
\end{equation}
From the second compatibility relation we get
\begin{equation}\label{4.2.4}
a_{n+1}^2T_{n+1}-a_n^2T_{n-1}=-b_n t_{n+1}+(b_n-1)t_n+N,
\end{equation}
\begin{equation}\label{4.2.5}
a_{n+1}^2-a_n^2=c(t_{n+1}-t_n-1).
\end{equation}
From (\ref{4.2.5}) we find, using telescopic summation, that
\begin{equation}\label{4.2.6}
a_n^2=c(t_n-n).
\end{equation}
Equation (\ref{4.2.3}) gives
\begin{equation}\label{4.2.7}
b_n=-1-c+n+\alpha+c T_n.
\end{equation}
After multiplying (\ref{4.2.4}) by $T_n$ we can replace $-b_nT_n$
using (\ref{4.2.2}), which allows us to take a telescopic sum. As a
result we get
\begin{equation}\label{4.2.8}
c^2(t_n-n)T_nT_{n-1}=c t_n^2+cN\sum_{j=0}^{n-1}T_j-c
t_n\sum_{j=0}^{n-1}T_j-t_n(cN-1+\alpha).
\end{equation}
One can also replace $a_n^2$ and $a_{n+1}^2$ in equation
(\ref{4.2.4}) using (\ref{4.2.6}) and rewrite the obtained equation
by  collecting the coefficients of $t_{n+1}$ and $t_n$. Taking a
telescopic sum and using (\ref{4.2.7}) gives
\begin{equation}\label{4.2.9}
c(t_n-n)(T_n+T_{n-1})+t_n(\alpha+n-c-2)=nN-c\sum_{j=0}^{n-1}T_j.
\end{equation}
Multiplying (\ref{4.2.9}) by $T_n$, using (\ref{4.2.8}) on the
left-hand side and (\ref{4.2.2}) on the right-hand side and
eliminating the sum using (\ref{4.2.2}) we finally obtain
\begin{multline}\label{4.2.11}
(N+c T_n)(1-\alpha-cN+c t_{n+1}+cT_n(\alpha-c-2+cT_n))\\+c
t_n(N-t_{n+1}+c T_n)=0.
\end{multline}
Replacing the sum in equation (\ref{4.2.9}) with the aid of
(\ref{4.2.8}) we get
\begin{multline}\label{4.2.12}
t_n(1-2n+(2n-2+ \alpha)(N+1)+c(n+N)T_n+c(n+N+cT_n)T_{n-1})\\= n(N+c
T_{n-1})(N+c T_n)+t_n^2(n-2+\alpha+c(T_n+T_{n-1})).
\end{multline}
Setting
\begin{equation}\label{subst}
x_n= \frac{t_n}{N},\quad y_n=- \frac{c T_n}{N}-1
\end{equation}
in the last two equations, we obtain the system \eqref{gen system}.

Finally, we note that the initial conditions for the recurrence
coefficients, in terms of the moments $\mu_1$ and $\mu_0$, are given
by $a_0^2=0 $ and
$$b_0 = \frac{ \mu_1 }{ \mu_0}
= \frac{c N}{(1-\alpha)}\frac{ M(-N+1, 2-\alpha, -c)} { M( -N,
1-\alpha, -c)},$$ where $M(a, b, z)$ is the confluent hypergeometric
function $_1F_1(a; b; z)$. This, together with \eqref{def:x_n} and
\eqref{def:y_n}, implies initial conditions \eqref{eq:initial value
xy} for $x_n$ and $y_n$.

This completes the proof of Theorem \ref{Th1}. 


\section{Proof of Theorem \ref{thm: P5}} \label{sec:proof of p5}


We use the method proposed in \cite{Meixner1, Meixner2} and system
(\ref{gen system}) to prove Theorem \ref{thm: P5}. We repeat the
main steps to be self-contained.

To prove the theorem, we derive the differential equation for $y_n$.
In \cite{Meixner1, Meixner2} we have used the Toda system
\begin{equation} \label{Toda}
\begin{cases} \ds \left(a_n^2\right)':=\frac{d}{dc}\left(a_n^2\right)=\frac{a_n^2}{c}(b_{n}-b_{n-1}), \\
\ds b_n':=\frac{d}{dc}b_n=\frac{1}{c}(a_{n+1}^2-a_n^2),
\end{cases}
\end{equation}
which also holds in the present case. Solving the first equation
(\ref{gen system}) for $x_{n+1}$ and the second equation for
$y_{n-1}$ and substituting into the Toda system (\ref{Toda}) (where
we have replaced $a_n^2$ and $b_n$ by their expressions in terms of
$x_n$ and $y_n$ from Theorem~\ref{Th1}), we get two equations
\begin{eqnarray*}x_n'  
&=&\frac{Nx_n-n}{c}\left(-x_n-y_n+\frac{x_n(-N-1+Nx_n)(\alpha-N-1+Nx_n)}{N(Nx_n-n)(x_n+y_n)}\right)
\end{eqnarray*}
 and
\begin{eqnarray}\label{for transofmration xn}
y_n'  
&=&\label{*}x_n+y_n+\frac{y_n(N+1+Ny_n)(N+1-\alpha+Ny_n)}{cN(x_n+y_n)},
\end{eqnarray}
where the differentiation is with respect to $c$. By differentiating
equation (\ref{*}) and substituting the expression for $x_n'$ we
obtain an equation for $y_n''$ as a function of $y_n',y_n,x_n$:
\begin{eqnarray*}y_n''&=&-\frac{1}{c}(Nx_n-n)(x_n+y_n)+\frac{x_n(-N-1+Nx_n)(\alpha-N-1+Nx_n)}{cN(x_n+y_n)}+y_n'\\
&+&\frac{y_n'}{Nc(x_n+y_n)}\left(
(N+1+Ny_n)(N+1-\alpha+Ny_n)+Ny_n(2N+2-\alpha+2Ny_n)\right)\\
&-&\frac{y_n}{Nc^2(x_n+y_n)^2}(N+1+Ny_n)(N+1-\alpha+Ny_n)\\
&\times&\left(x_n+y_n+cy_n'-(Nx_n-n)(x_n+y_n)
+\frac{x_n(-N-1+Nx_n)(\alpha-N-1+Nx_n)}{N(x_n+y_n)}\right).\end{eqnarray*}

Eliminating $x_n$ between this equation and (\ref{*}) gives a
nonlinear second order second degree equation for $y_n$:
\begin{equation}\label{nonlinear diff vn}
G(y_n'',y_n',y_n,c)=0.
\end{equation}
We have used \texttt{Mathematica}\footnote{www.wolfram.com} to
compute this long expression. It can be checked by direct
computations, that applying a transformation
$$y_n(c)=\frac{v(z)+2\alpha-4N+2n-3}{4N}, \;\;c=z^2$$ this equation
simplifies considerably and can be written as
$$(v''-6v^2-\alpha_1 v-\beta_1)^2=(v/z-2z)^2(v'^2-4v^3-\alpha_1
v^2-2\beta_1 v-\gamma_1),$$ where $v=v(z)$ and
\begin{gather*}
\alpha_1=4(2\alpha-4N+6n-1),\\
\beta_1=2(2n+1)(6n-8N-5)+8(2n+1)\alpha-8\alpha^2,\\
\gamma_1=4(2\alpha-4N+2n-3)(4n^2+4n+1-4\alpha^2).
\end{gather*}
 The last equation appears in \cite{Cosgrove} (equation (A.8)) and is known to be related to the fifth Painlev\'e equation. This is similar to the case of generalized Meixner polynomials (see \cite{Meixner1}).

By taking $y_n(c)$ in a form as shown in \eqref{transformation}, we
get the fifth Painlev\'e equation \pf \ \eqref{P5} with parameters
\eqref{par}.

A similar approach can be used to obtain a (second order second
degree) differential equation for
$x_n$. 
To get the expression \eqref{transformation xn} for $x_n(c)$ in terms of the solutions of the fifth Painlev\'e equation we use \eqref{for transofmration xn}, substitute \eqref{transformation} and find $x_n(c)$ by solving a quadratic equation. One can check with \texttt{Mathematica} that one of the roots of this equation, i.e. \eqref{transformation xn}, indeed gives the statement of Theorem \ref{thm: P5}. This completes the proof of Theorem \ref{thm: P5}. 

Note that we can  derive a nonlinear discrete second order equation
for $y_n(c)$. From the first equation of (\ref{gen system}) with $n$
and the second equation with $n+1$ we eliminate $x_{n+1}$ by
computing the resultant. The obtained equation and the second
equation of (\ref{gen system}) with $n$ then allow us to eliminate
$x_n$. As a result, we obtain a nonlinear discrete equation for
$y_n(c)$ which we denote by
\begin{equation}\label{nonlinear discrete vn}
F_1(y_{n-1},y_{n},y_{n+1},c)=0.
\end{equation}
The equation was again obtained by using \texttt{Mathematica} but it
is too long and too complicated to include here explicitly.
Similarly, a nonlinear discrete equation for $x_n(c)$ can be
obtained: \begin{equation*}
F_2(x_{n-1},x_n,x_{n+1},c)=0.\end{equation*} One can, similarly to
\cite[Sect. 3]{Meixner1}, check that the functions $y_{n-1}$,
$y_{n}$ and $y_{n+1}$ in \eqref{nonlinear discrete vn} are connected
by using the B\"acklund transformation of the fifth Painlev\'e
equation \cite[$\S$ 32.7 (v)]{DLMF}. In particular, one can express
all of them using only $y$ and $y'$ (the solution of \pf\ with
parameters \eqref{par}). When these expressions are substituted into
\eqref{nonlinear discrete vn}, it becomes identically zero.
Moreover, equation \eqref{nonlinear discrete vn} can  essentially
(up to a factor depending on  $y_n$) be obtained by eliminating $y$
between
$$
y_{n+1}=\frac{-(1+N+Ny y_n)(N-n+(1+N-\alpha+N y_n)y)}{N(N-n+(1+n+N
y_n)y)}
$$
and
$$
y_{n-1}= \frac{(1+c+N-\alpha+N y_n-(1+N-\alpha+N y_n)y)P}{N(y-1)( n
c+N (1+c+N-\alpha+N y_n-(1+N-\alpha+N y_n)y)y_n)},
$$
where $P=(N+1)c+N (c-n+2+2N-\alpha+N y_n+(n-2N-2+\alpha-N
y_n)y)y_n.$

\begin{figure}[h]
\centering \resizebox{2in}{!}{\includegraphics{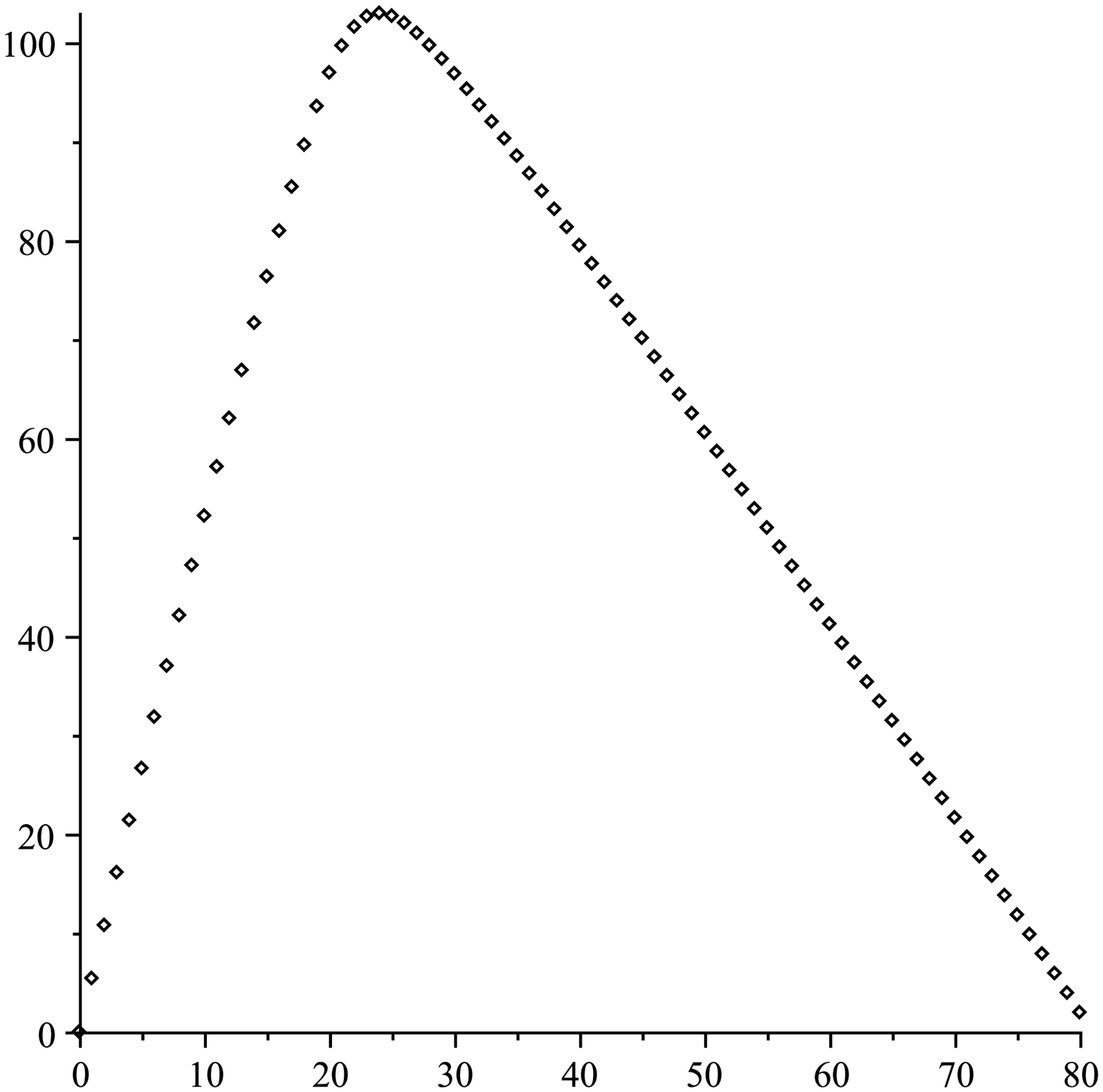}}
\resizebox{2in}{!}{\includegraphics{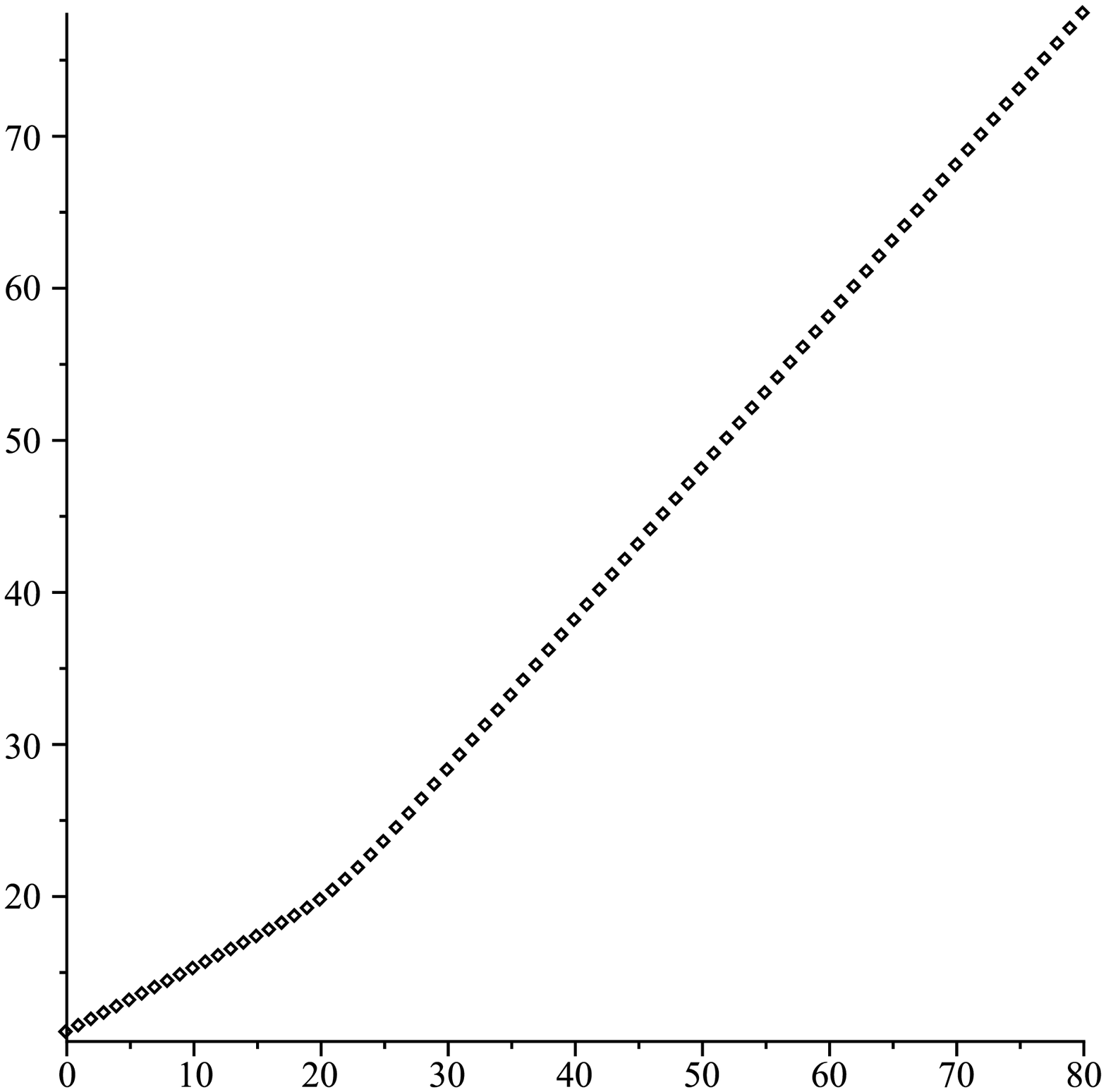}}
\caption{Recurrence coefficients ($a_n^2$ on the left, $b_n$ on the
right) for generalized Krawtchouk polynomials ($N=80$, $\alpha=-1$,
$c=2$)} \label{fig:1}
\end{figure}


\begin{figure}[h]
\begin{center}
\subfigure[]{
\resizebox*{5cm}{!}{\includegraphics{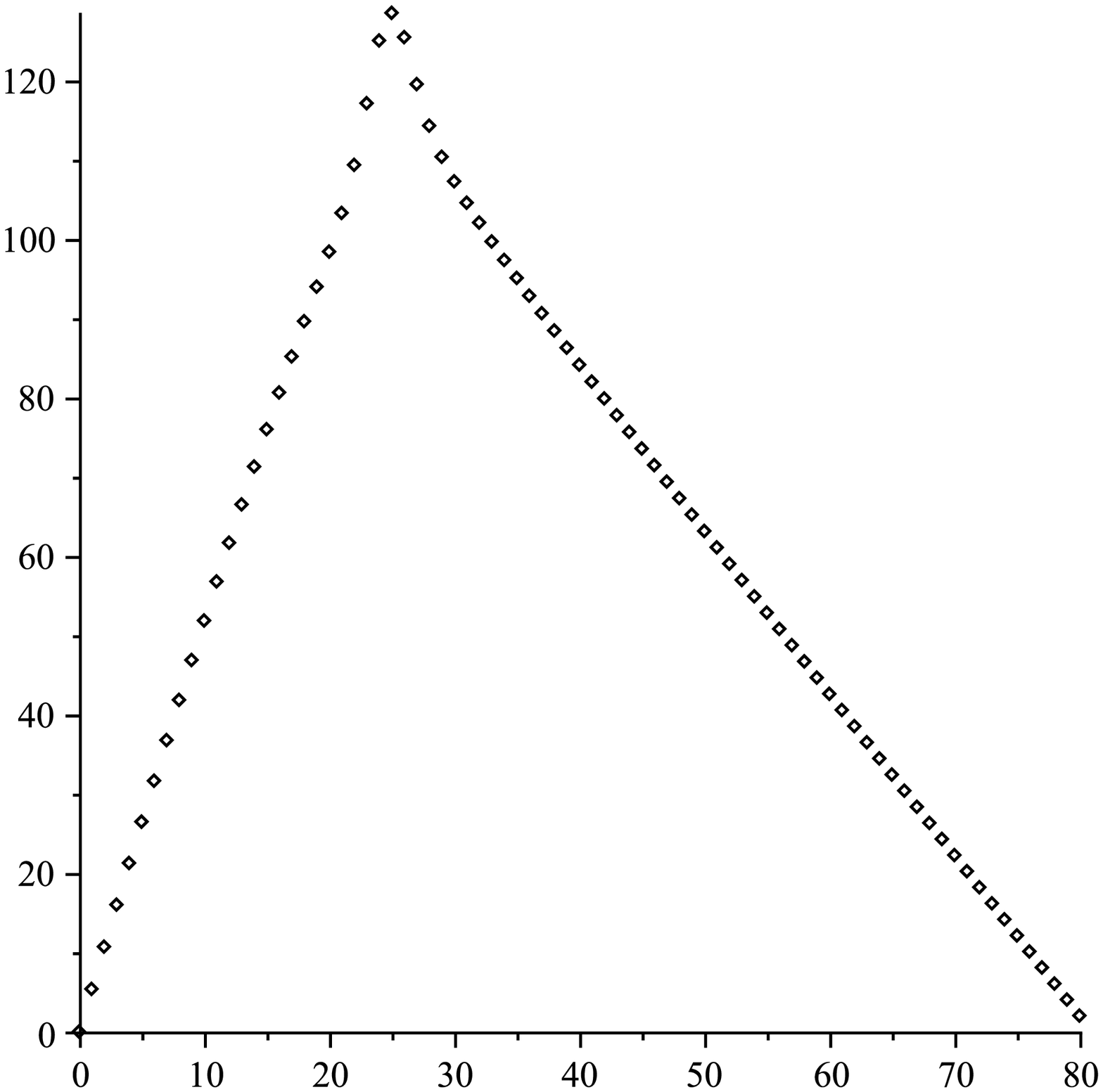
}}} \subfigure[]{
\resizebox*{5cm}{!}{\includegraphics{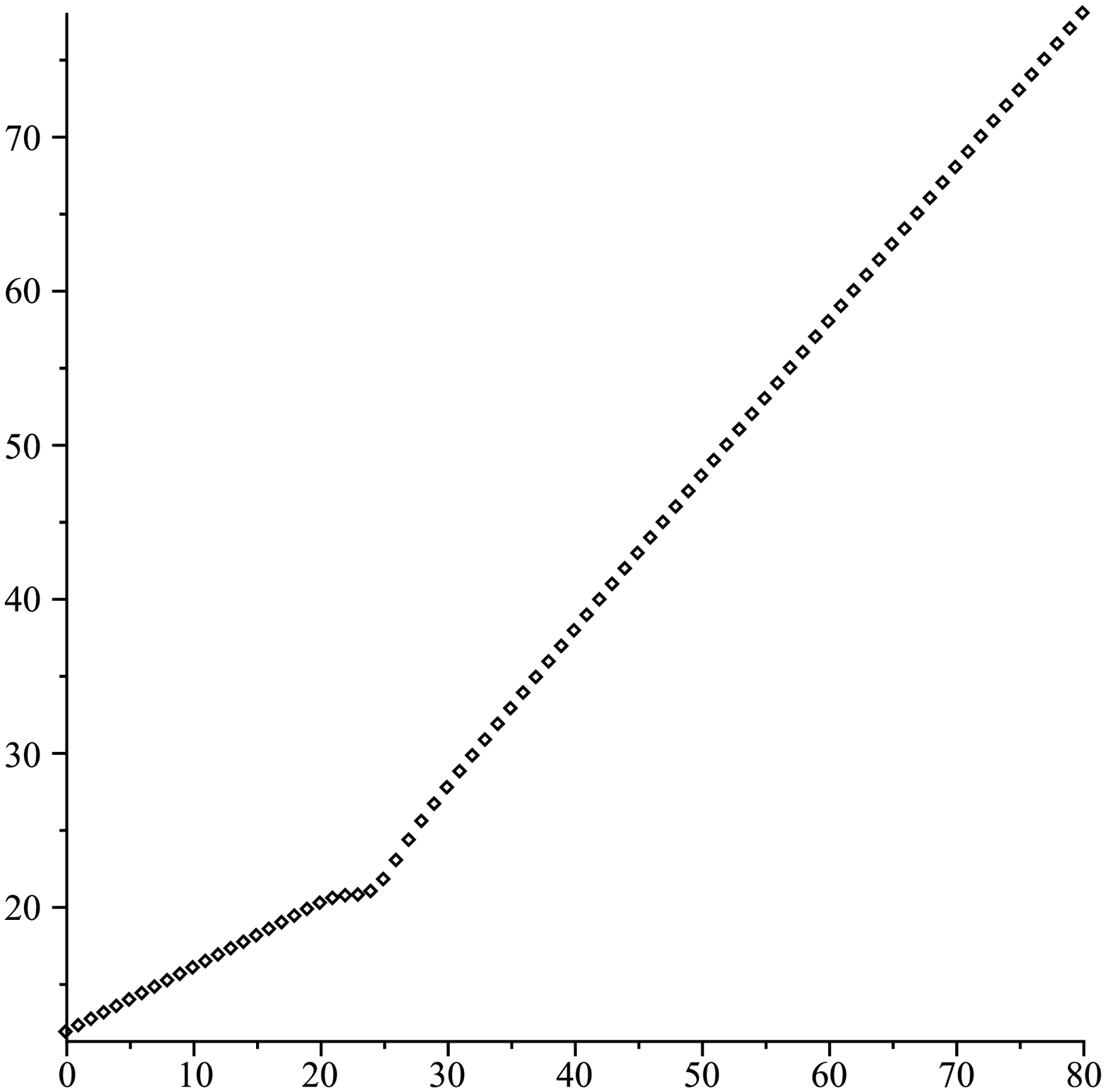
}}} \caption{Recurrence coefficients ($a_n^2$ on the left, $b_n$ on
the right) for generalized Krawtchouk polynomials ($N=80$,
$\alpha=0.8$, $c=2$).} \label{fig:2}
\end{center}\end{figure}

\section{Computing the recurrence coefficients}\label{num}


\begin{figure}[h]
\begin{center}
\subfigure[]{
\resizebox*{5cm}{!}{\includegraphics{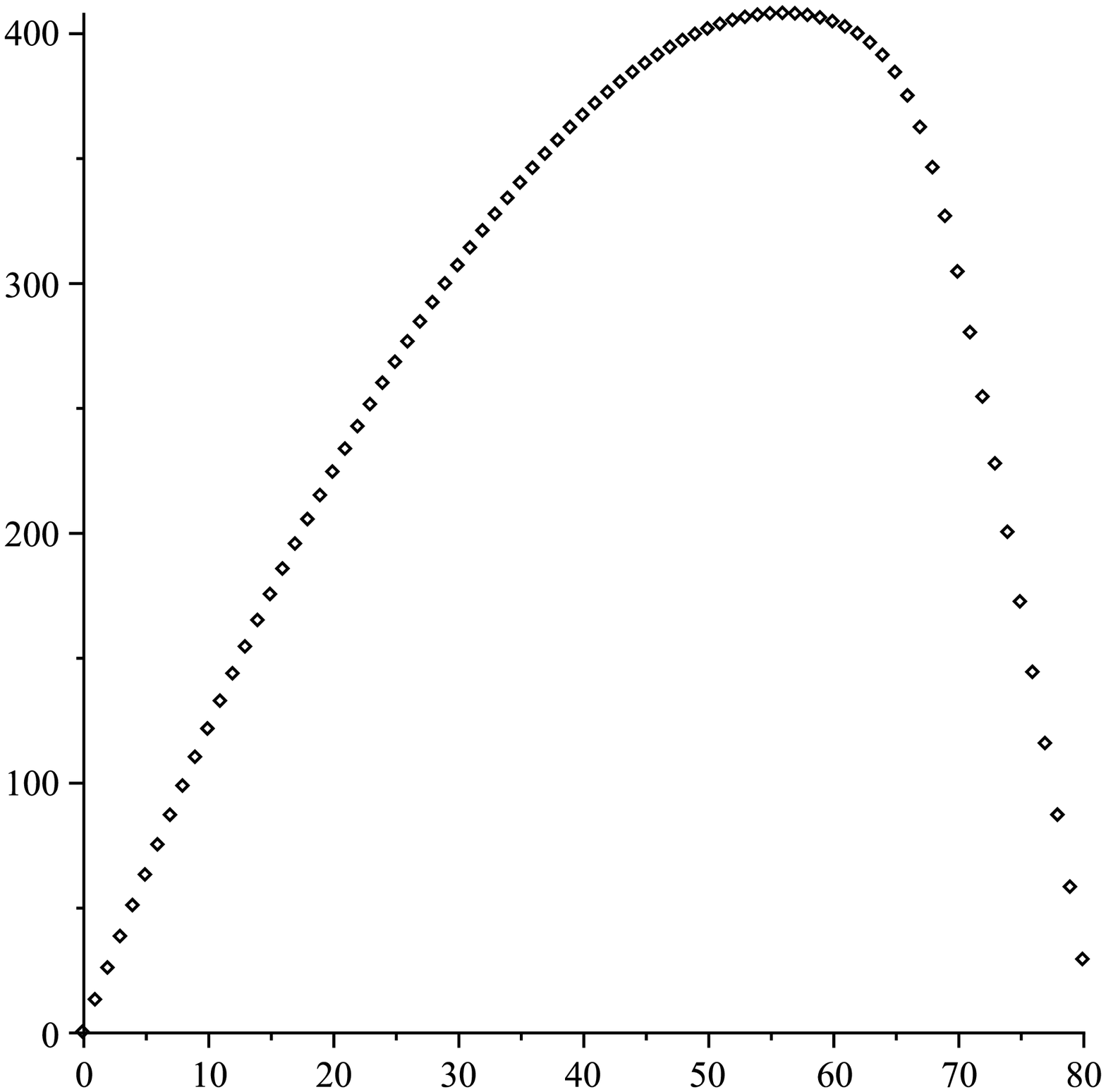
}}} \subfigure[]{
\resizebox*{5cm}{!}{\includegraphics{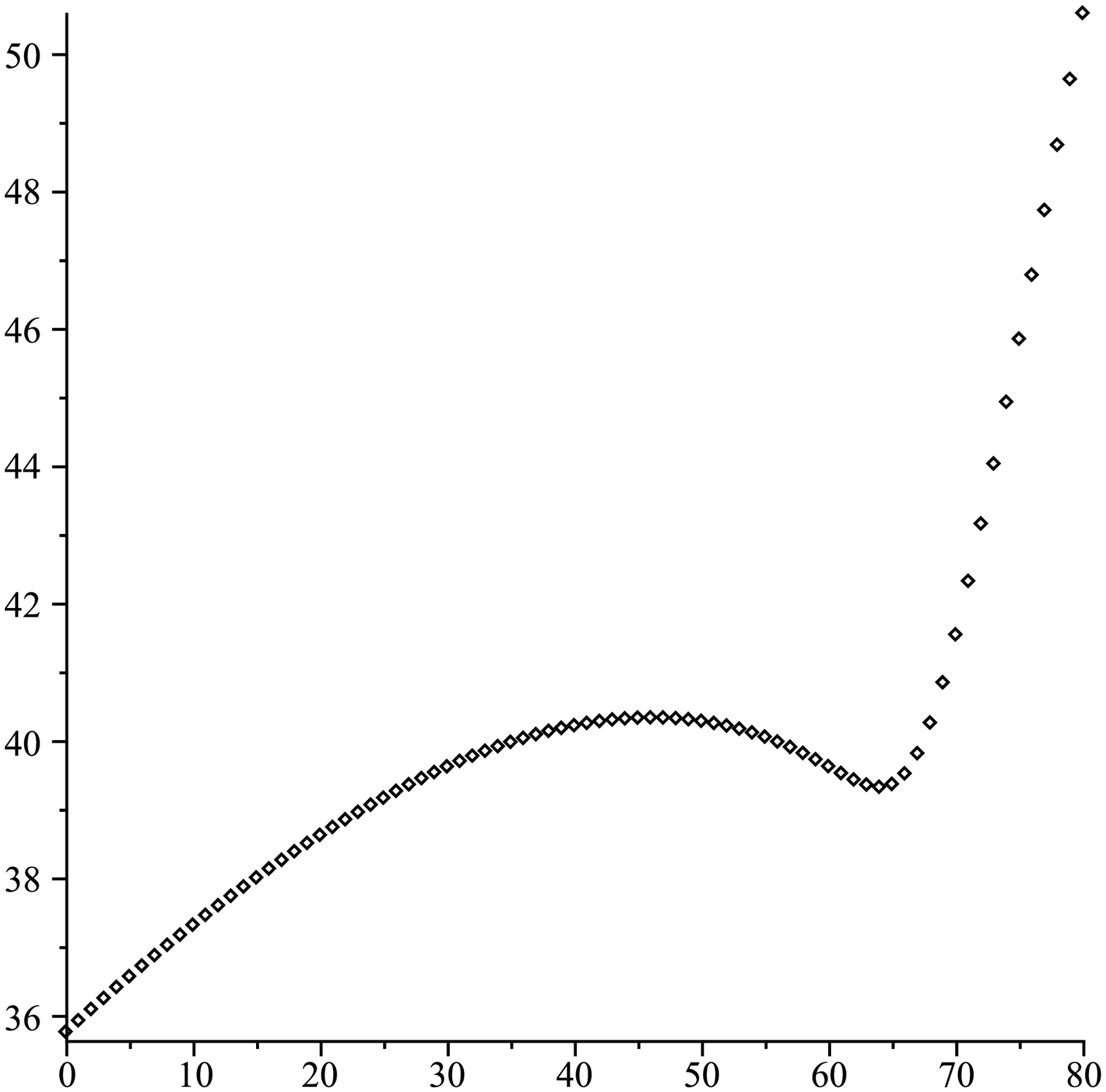
}}} \caption{Recurrence coefficients ($a_n^2$ on the left, $b_n$ on
the right) for generalized Krawtchouk polynomials ($N=80$,
$\alpha=-1$, $c=30$).} \label{fig:3}
\end{center}\end{figure}

\begin{figure}[h]
\begin{center}
\subfigure[]{ \resizebox*{5cm}{!}{\includegraphics{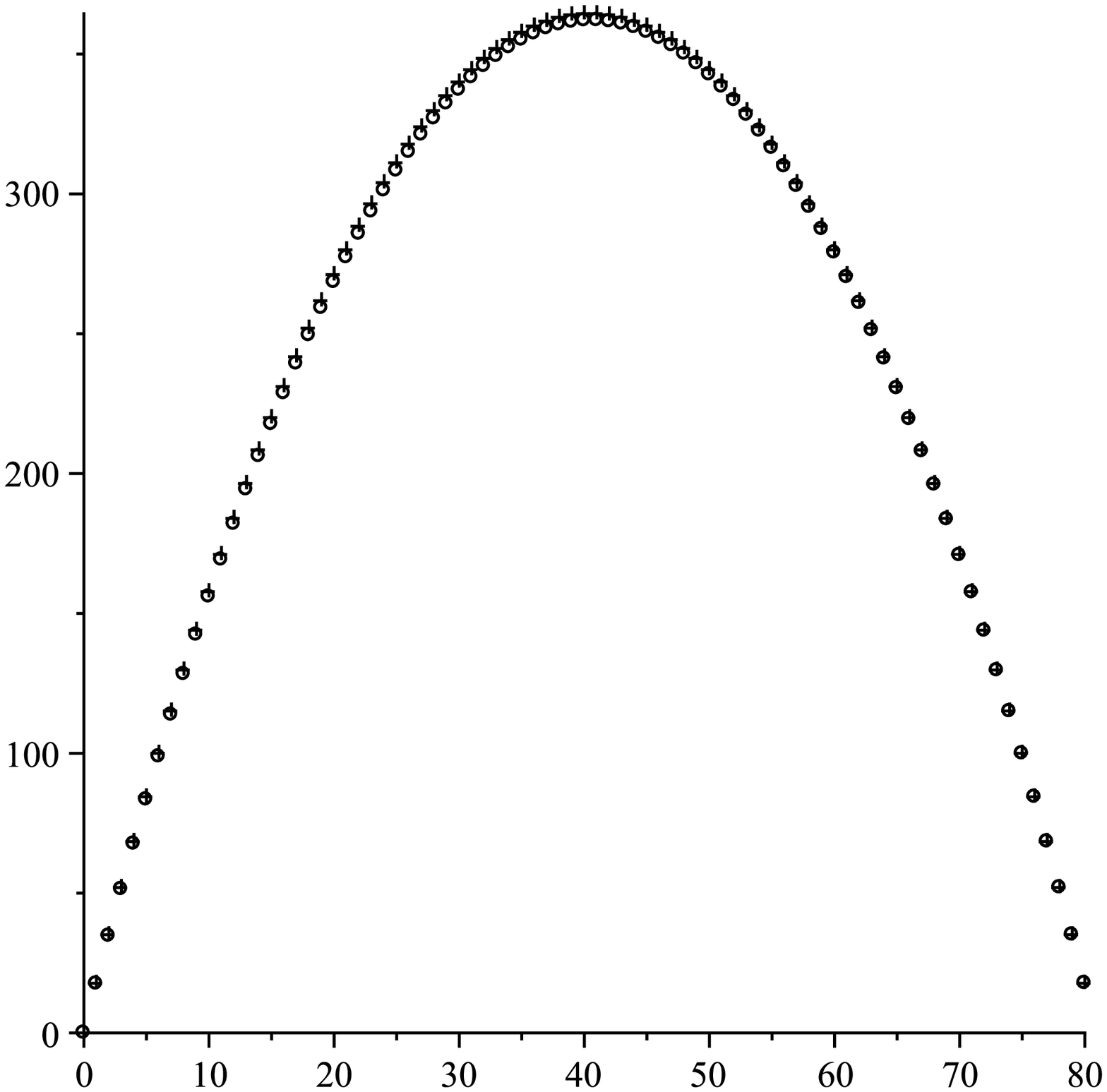}}}
\subfigure[]{ \resizebox*{5cm}{!}{\includegraphics{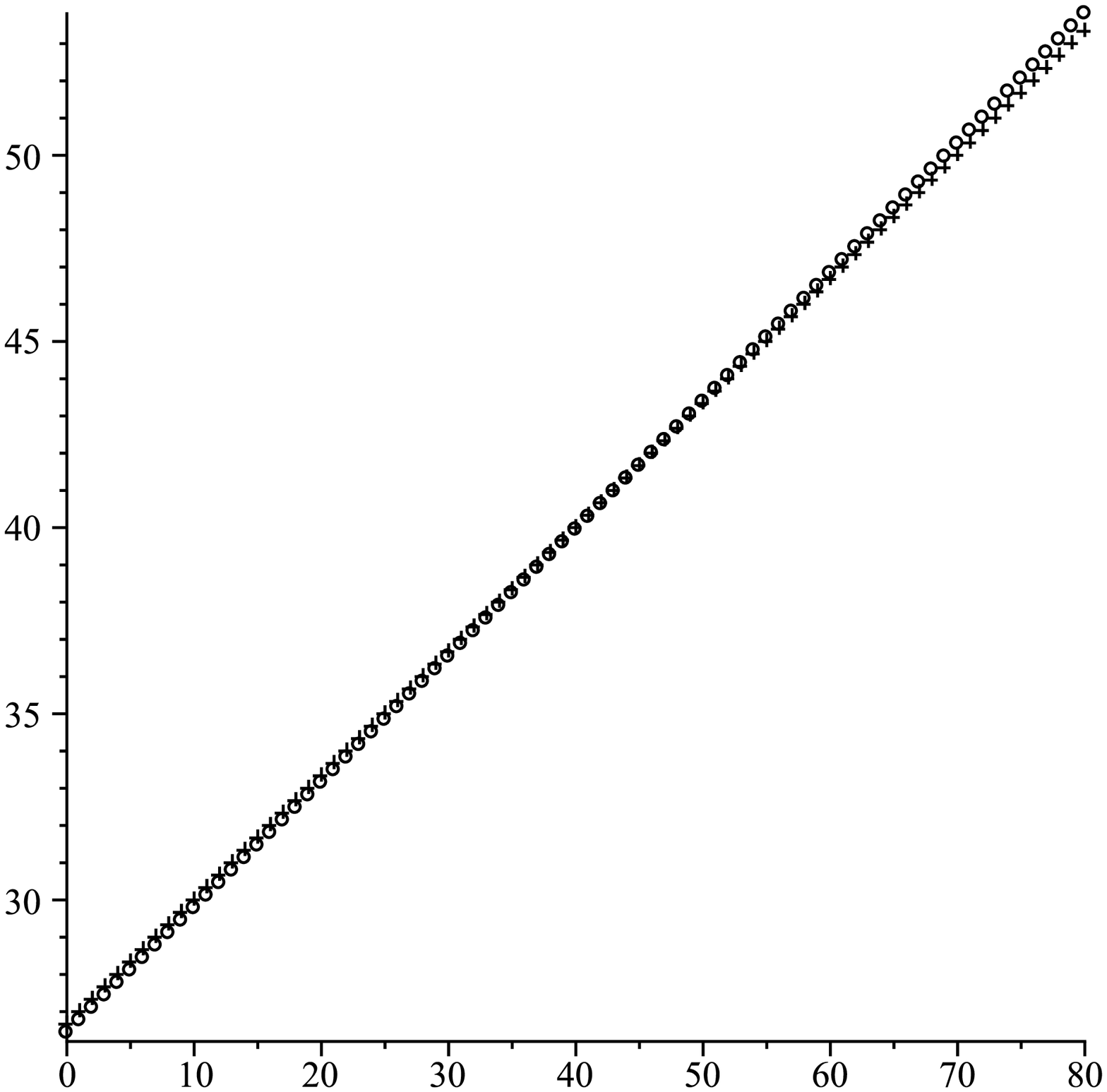}}}
\caption{Recurrence coefficients ($a_n^2$ on the left, $b_n$ on the
right) for generalized Krawtchouk polynomials ($N=80$,
$\alpha=-2000$, $c=1000$, circles) and for Krawtchouk polynomials
($N=80$, $p=1/3$, crosses)} \label{fig:4}
\end{center}\end{figure}

The system (\ref{gen system}), together with the definitions
(\ref{def:x_n}) and (\ref{def:y_n}), can be used to compute the
recurrence coefficients starting with the initial values $x_0$ and
$y_0$ given in (\ref{eq:initial value xy}).
  Figures 1-3 show these coefficients for
$N=80$ and for some different choices of the parameters $c > 0$ and
$\alpha < 1$.
Recall that the recurrence coefficients (\ref{KrawRec}) for the
Krawtchouk polynomials
 are such that $a_n^2$ is quadratic in $n$, with $a_0^2=0=a_{N+1}^2$,
and that $b_n$ is of degree 1 in $n$. Also recall that when
 $c$ and $-\alpha$ tend to infinity and $-\frac{c}{\alpha}\to\frac{p}{1-p}$
we recover the Krawtchouk weight. Figure~\ref{fig:4} shows this
observation.  For the generalized Krawtchouk polynomials the $a_n^2$
are first increasing until they reach their maximum, and then
decreasing, with $a_0^2=0=a_{N+1}^2$. The $b_n$ are no longer linear
but change from concave to convex. The point where $a_n^2$ reaches
its maximum and the inflection point for the $b_n$ depend on the
parameter $c$, these two points do not coincide.  Moreover, for
$\alpha$ values close to $1$, the $b_n$ change once more to concave,
as can be seen (albeit not very clearly) in Figure~\ref{fig:2}.

To get these values (they were obtained using \texttt{Maple}), one
has to work with quite high accuracy, for it turns out that even a
very small perturbation ($10^{-100}$ for $N=80$) in the initial
value $y_0$ quickly leads to very bad results.  An explanation for
this behaviour is the following: the system (\ref{gen system})
allows to compute the $x_n$, from which the coefficients $a_n^2$ can
be obtained by (\ref{def:x_n}).  There is, however, no a priori
reason why these $a_n^2$ should be positive and $a_{N+1}^2=0$.  This
is not a new observation: in \cite{LiesThesis} (Section 6.5, about
discrete $q$-Hermite I polynomials) it was conjectured that the
initial value needed to obtain the recurrence coefficients for these
polynomials, is the only real value which leads to positive values
for $a_n^2$ for every $n>0$. Since in this generalized Krawtchouk
case there is only a finite number of coefficients
($a_1^2,\ldots,a_N^2$) which need to be positive, this conjecture
obviously does not hold here, but we conjecture that there is only
one initial value $y_0$ for which $a_{N+1}^2=0$ and the
$a_1^2,\ldots,a_N^2$ are positive. An alternative, and numerically
more stable way, would be to use a fixed point algorithm for which
the fixed point satisfies (\ref{gen system}).

Following a referee's remark, we also used the Stieltjes approach to
calculate the recurrence coefficients, see e.g. \cite[\S
2.2.3.1]{Gautschi}. In this method, one uses
(\ref{recu:orthonormal}) to calculate $a_n$ and $b_n$ from $p_n$ and
$p_{n-1}$.  Then (\ref{3 term orthonormal}) can be used to obtain
$a_{n+1}p_{n+1}(x)$, which after normalization gives $p_{n+1}$.  The
initial values needed are $p_{-1}=0$ and
$p_0=\frac{1}{\sqrt{\sum_{k=0}^Nw(k)}}$. Numerical calculations show
that this method needs less precision in the initial conditions but
requires more computations.


\section{Conclusion and future directions}\label{sec:dis}

In this paper we have been dealing with the recurrence coefficients
of the discrete orthogonal polynomials, namely the generalized
Krawtchouk polynomials, and have shown that they are related to the
classical solutions (in terms of confluent hypergeometric functions)
of  the fifth Painlev\'e equation. A possible further direction is
to consider more factors in \eqref{eq:potential} and find out
whether the sixth Painlev\'e equation or some equation from the
Painlev\'e hierarchy is related to the recurrence coefficients of
such a weight. A reason that leads to this conjecture follows from
the observation that the initial conditions for such a weight
correspond to special solutions of the sixth Painlev\'e equation.

\section*{Acknowledgements}
The authors are grateful to the referees for providing a lot of
useful and helpful suggestions which substantially improved the
presentation of this paper.

This paper was started while GF was visiting KU Leuven in June 2011.
The hospitality of the Department of Mathematics is gratefully
acknowledged. GF is also supported by the Polish MNiSzW Iuventus
Plus grant Nr 0124/IP3/2011/71 and is partially supported by MNiSzW
Grant N N201 397937. This work was supported by FWO project
G.0427.09 and KU Leuven research grant OT/08/033. LZ is a
Postdoctoral Fellow of the Scientific Research Foundation - Flanders
(FWO), Belgium. The authors are grateful to Peter Clarkson for
illuminating discussions.

\begin{verbatim}
Galina Filipuk
Faculty of Mathematics, Informatics and Mechanics
University of Warsaw
Banacha 2 Warsaw 02-097, Poland
filipuk@mimuw.edu.pl


Lies Boelen, Christophe Smet, Walter Van Assche, Lun Zhang
Department of Mathematics
KU Leuven
Celestijnenlaan 200 B box 2400
BE-3001 Leuven
Belgium
lies.boelen@gmail.com
christophe@wis.kuleuven.be
walter@wis.kuleuven.be
lun.zhang@wis.kuleuven.be
\end{verbatim}

\end{document}